\def\FormatStyle{InHouse}  
\def\InHouse{InHouse}
\def\Conference{Conference}

\newif\ifblind
\blindfalse 

\ifx \FormatStyle \InHouse  
  \documentclass[a4paper,12pt,onecolumn]{article}
  \usepackage[margin=1in]{geometry}
  \usepackage{setspace}
  \usepackage{fancyhdr}
\else 
  \ifx\FormatStyle\Conference  
  \documentclass[a4paper,10pt,twocolumn]{article}
  \else 
  \documentclass[a4paper,10pt,onecolumn]{article}
  \fi
  \usepackage[margin=1in]{geometry}
\fi

\usepackage[utf8]{inputenc}
\usepackage[noblocks,affil-it]{authblk}
\usepackage[cmex10]{amsmath}
\usepackage{amsthm,amsfonts,amssymb}
\usepackage[verbose]{wrapfig}
\usepackage{IEEEtrantools}  
\usepackage[pdftex]{graphicx}
\usepackage[usenames,dvipsnames,svgnames,table]{xcolor}
\graphicspath{{./Figures/}}
\usepackage{appendix}
\usepackage{xspace}

\usepackage{algorithm,algpseudocode}
\algrenewcommand{\algorithmiccomment}[1]{\hskip3em - #1}

\ifx\FormatStyle\InHouse
 \usepackage[disable]{todonotes}
\else
 \usepackage[disable]{todonotes}
\fi

\usepackage{natbib} 
\usepackage{hyperref} 
\definecolor{Pcolor}{rgb}{0.20,0.50,0.40}
\hypersetup{colorlinks=true, citecolor = Pcolor}

\newtheorem{defn}{Definition}

\newtheorem{assumption}{Assumption}

\newcommand{\mynotesH}[2][]{\todo[backgroundcolor=blue!20!white,inline,
bordercolor=red]{#2}}

\newcommand{\bbR}{\mathbb{R}}

\newcommand{\I}{\mathbf{I}}
\newcommand{\Hm}{\mathbf{H}}

\newcommand{\X}{\mathbf{X}}

\newcommand{\Z}{\mathbf{Z}}

\newcommand{\Rs}{R$^2$\xspace}
\newcommand{\mRs}{\text{\Rs}} 

\newcommand{\numspace}{\IEEEeqnarraynumspace}

\DeclareMathOperator*{\argmax}{arg\,max}
\newcommand\addtag{\refstepcounter{equation}\tag{\theequation}}

\ifblind
  \author{}
\else
  \author{Kory D. Johnson}
  \author{Robert A. Stine}
  \author{Dean P. Foster}
  \affil{The Wharton School\\University of Pennsylvania}
\fi

\ifx \FormatStyle \InHouse
 \title{Submodularity in Statistics: Comparing the Success of Model
Selection Methods}
 \date{\today}
\else
 \title{Submodularity in Statistics: Comparing the Success of Model Selection
Methods}
 \date{}
\fi

\begin{document}

\maketitle
\listoftodos

\begin{abstract}

We demonstrate the usefulness of submodularity in statistics as a
characterization of the difficulty of the \emph{search} problem of feature
selection. The search problem is the ability of a procedure to identify an
informative set of features as opposed to the performance of the optimal
set of features. Submodularity arises naturally in
this setting due to its connection to combinatorial optimization. 
In statistics, submodularity isolates cases in which collinearity makes the
choice of model features difficult from those in which this task is routine.
Researchers often report the signal-to-noise ratio to measure the
difficulty of simulated data examples. A measure of submodularity should also be
provided as it characterizes an independent component difficulty. Furthermore,
it is closely related to other statistical assumptions used in the development
of the Lasso, Dantzig selector, and sure information screening.



\end{abstract}


\section{Introduction}
\label{sec:intro}

We study the problem of selecting predictive features from a large feature
space. Our data consists of $n$
observations of (response, feature) sets, ($y_i,x_{i1},\ldots,x_{im}$), where
each observation has $m$ associated features. Observations are collected into
matrices and the following model is assumed for our data
\begin{IEEEeqnarray}{rCl+rCl}
 Y & = & \X\beta + \epsilon & \epsilon & \sim & N_n(\mathbf{0},\sigma^2 \I_n)
 \label{eqn:model}
\end{IEEEeqnarray}
where $\X$ is an $n \times m$ matrix and $Y$ is an $n \times 1$ response vector.
Typically, most of the elements of $\beta$ are 0. Hence, generating good
predictions requires identifying the small subset of predictive features. The
model (\ref{eqn:model}) proliferates the statistics and machine learning
literature. In modern applications, $m$ is often large, potentially with $m \gg
n$, which makes the selection of an appropriate subset of these features
essential for prediction.

The model selection problem is to minimize the error sum of squares
\[\text{ESS}(\hat Y) = \|Y - \hat Y\|^2_2 = \sum_{i=1}^n (Y_i - \hat Y_i)^2\]
while restricting the number of nonzero coefficients:
\begin{IEEEeqnarray}{c't'rCl}
 \min_\beta \text{ESS}(\X\beta) & s.t. &
\|\beta\|_{l_0} = \sum_{i=1}^m I_{\{\beta_{i}\neq0\}} \leq k,
\label{eqn:sparse-reg}
\end{IEEEeqnarray}
where the number of nonzero coefficients, $k$, is the desired sparsity. Note
that we are not assuming a sparse representation exists, merely asking for a
sparse approximation. In the statistics literature, the model selection
problem (\ref{eqn:sparse-reg}) is more commonly posed as a penalized regression:
\begin{IEEEeqnarray}{rCl}
  \label{eqn:penreg}
\hat\beta_{0,\lambda} & = & \text{argmin}_{\beta}\left\{\text{ESS}(\X \beta) +
  \lambda\|\beta\|_{l_0} \right\}
\end{IEEEeqnarray}
where $\lambda \ge 0$ is a constant. The classical hard thresholding algorithms
$C_p$ \citep{Mal73}, AIC \citep{Aka74}, BIC \citep{Schwarz78}, and RIC
\citep{FosterG94} vary $\lambda$. The solution to (\ref{eqn:penreg})
is the least-squares estimator on an optimal subset of features. Let $M
\subset \{1,\ldots,m\}$ indicate the coordinates of a given model so that
$\X_M$ is the corresponding submatrix of the data. If $M_\lambda^*$ is the
optimal set of features for a given $\lambda$ then $\hat\beta_{0,\lambda}^* =
(\X_{M^*}^T \X_{M^*})^{-1}\X_{M^*}Y$.

Given the combinatorial nature of the constraint, solving (\ref{eqn:sparse-reg})
quickly becomes infeasible as $m$ increases and is NP-hard in general
\citep{Nata95}. Forward stepwise is the greedy approximation to the
solution of (\ref{eqn:sparse-reg}). Let $M_i$ be the features in the forward
stepwise model after step $i$ and note that the size of the model is $|M_i| =
i$. The algorithm is initialized with $M_0 = \emptyset$ and iteratively adds
the variable which yields the  largest reduction in ESS. Hence, 
$M_{i+1} = \{M_i\cup j\}$ where
\[j = \argmax_{l \in \{1,\ldots,m\}\backslash M_i}
\text{ESS}(\X_{M_i\cup l} \hat\beta_{M_i\cup l}^\text{LS}).\]
After the first feature is selected, subsequent models are built having fixed
that feature in the model. $M_1$ is the optimal size-1 model, but $M_i$ for
$i\geq 2$ is not guaranteed to be optimal, because $M_i$ is forced to include
the features identified at previous steps.

Subset selection problems are difficult because features can interact in
unexpected ways. Here,
``unexpected'' means that the change in model fit when adding a feature can be
completely different depending on the other features in the model. This paper
characterizes the cases in which features produce such unexpected results. 


A simple example from \citet{Miller02} clarifies this point. Suppose
forward stepwise is run on the data in Table
\ref{tab:simple-ex}. The first step selects the feature that
is maximally correlated with $Y$. For features $\X_1$, $\X_2$, and $\X_3$,
these are $r_{Y1} = .0$, $r_{Y2} = -.0016$, and $r_{Y3} = .4472$, respectively.
Therefore, forward stepwise selects $\X_3$ on the first step.
The second step chooses the feature with the maximum partial correlation. That
is, the maximum correlation when features are considered orthogonally to
$\X_3$. For $\X_1$ and $\X_2$ these are $r_{Y1.3^\perp} = .0$ and
$r_{Y2.3^\perp} = -.0014$, respectively. Forward stepwise appears to find a
significant features on the first step, but then no other features  seem
important. The true equation for the response, however, is $Y = \X_1 - \X_2$.
This cannot be identified by forward stepwise because of the ``incorrect'' first
step which includes $\X_3$. Furthermore, unless forward stepwise
continues to select features even when they appear uninformative,
the optimal set can never be found. Intuitively, the difficulty arises because
$\X_1$ and $\X_2$ have large errors which cancel out.
\begin{table}[h]
\centering
 \caption{Simple data in which forward stepwise fails to identify the correct
model.}
 \label{tab:simple-ex}
  \begin{tabular}{cccc}
    $Y$ &$\X_1$ & $\X_2$ & $X_3$\\\hline
    -2 & 1000 & 1002 & 0 \\
    -1 & -1000 & -999 & -1 \\
    1 & -1000 & -1001 & 1 \\
    2 & 1000 & 998 & 0
  \end{tabular}
\end{table}

Most of our discussion concerns maximizing the model fit as opposed to
minimizing loss. Let $[m] = \{1,\ldots,m\}$. For a subset of indices $S \subset
[m]$, we denote the corresponding columns of our data matrix as $\X_S$, or
merely $S$ when the overloaded notation will not cause confusion. Our measure of
model fit for a set of features $\X_S$ is the coefficient of determination, \Rs,
defined as \[\mRs(S) = 1 - \frac{\text{ESS}(\X_S \hat\beta_S)}{\text{ESS}(\bar
Y)} \] where
$\bar Y$ is the constant vector of the mean response and $\hat\beta_S$ is the
least squares estimate of $\beta_S$. 

Forward stepwise performs well when the improvement in fit obtained by adding
a set of features to a model is upper bounded by the sum of the improvements
of adding the features individually. If a set of features improves the model
fit when considered together, a subset of those features must improve the fit
as well. Consider the improvement in fit by adding $\X_S$ to the model $\X_T$:
\[\Delta_T(S) := \mRs(S \cup T) - \mRs(T).\]Letting $S = A \cup B$, bound
$\Delta_T(S)$ as
\begin{IEEEeqnarray}{rCl}
\Delta_T(A) + \Delta_T(B) & \geq & \Delta_T(S). \label{eqn:submod.intuition3} 
\end{IEEEeqnarray}

If $A\cup B$ improves the model fit, equation (\ref{eqn:submod.intuition3})
requires that either $A$ or $B$ improve the fit when considered in isolation.
Therefore, signal that is present due to complex relationships among features
cannot be completely hidden when considering subsets of these features.
Equation (\ref{eqn:submod.intuition3}) defines a submodular function:
\begin{defn}[Submodular Function]
\label{def:submod-1}
  Let $F: 2^{[m]} \rightarrow \bbR$ be a set function defined on the the power
set of $[m]$. $F$ is submodular if $\forall A,B \subset [m]$
  \begin{IEEEeqnarray}{rCl}
    F(A) + F(B) & \geq & F(A \cup B) +F(A \cap B) \label{eqn:submod3}
  \end{IEEEeqnarray}
\end{defn}
The intuition in equation (\ref{eqn:submod.intuition3}) is recovered
by considering $A \cap B = \emptyset$. Alternatively, Definition
\ref{def:submod-1} can be rewritten as 
\begin{IEEEeqnarray*}{rCl}
 F(A) - F(A\cap B) + F(B) - F(A \cap B) & \geq & F(A \cup B) - F(A \cap B)\\
 \Rightarrow \Delta_{A\cap B}(A) + \Delta_{A\cap B}(B) 
  & \geq & \Delta_{A\cap B}(A\cup B),
\end{IEEEeqnarray*}
which considers the impact of $A \backslash B$ and $B \backslash A$ given $A
\cap B$. The influence of the union of set differences is less than the impact
of the sum of their marginal influences. We will refer to data as being
submodular if \Rs is a submodular function on the data.

Forward selection is a natural algorithm under submodularity as it adds the
feature to the model that yields the maximum marginal increase in fit. To fix
notation, if $S_i$ is the model at step $i$, feature $X_j$ is added if 
\[j = \arg\!\max_{l\notin S_i} \Delta_{S_i}(X_l)\]
and $S_{i+1} = \{S_i \cup j\}$.
Other greedy procedures have been proposed that change the criteria being
maximized at each step. For different criteria, this yields orthogonal matching
or orthogonal projection pursuit \citep{Barron+08,Miller02}. 

For all such methods, let $\hat Y^{(k)} = \X_{S_k}\hat\beta_{S_k}$ be the
estimated response after $k$
steps of the algorithm. Previous analyses determined the
rate at which ESS$(\hat Y^{(k)})$ decreases as a function of $k$
\citep{Barron+08,Jaggi13}. Instead, we focus on identifying the data conditions
that guarantee that ESS$(\hat Y^{(k)})$ is close to that of the optimal size $k$
subset. If forward stepwise is used and \Rs is submodular,
the classic result of \citet{NemWF78} shows that $\mRs(S_k)
\geq (1 - 1/e) \mRs(S_k^*)$, where $S_k^*$ is the subset of features which
solves the sparse regression problem in equation (\ref{eqn:sparse-reg}).

Instead of asking for an approximate solution to (\ref{eqn:penreg}), one can
relax the problem formulation. For example, the $l_0$ penalty can be
relaxed to an $l_1$ penalty, yielding the Lasso \citep{Tib96}.
Additionally, loss be measuring with the $l_\infty$ norm, which yields the
Dantzig selector \citep{CanT07}. While subset selection and greedy methods like
forward stepwise are classically studied, these relaxations have been the
primary focus of research in recent years. For cases where $\log(p) = O(n^{c})$
for $c>0$, the computational improvements from relaxing the constraint in
equation (\ref{eqn:sparse-reg}) do not produce efficient algorithms. In these
cases, a feature screening method can be used to reduce the dimensionality $p$
to feasible ranges before performing model selection \citep{FanL08}. 

These define two classes of algorithms: the first maintains the problem
formulation in (\ref{eqn:sparse-reg}) and provides approximate solutions, while
the provides exact solutions to relaxed problem formulations. Given
that both classes of algorithms can be used to answer the same question,
it is natural to ask which style of approximation is preferred. A general
framework comparing these as penalized regressions is given in \citet{FanL01},
and cases in which approximating (\ref{eqn:sparse-reg}) is preferable  to
solving the relaxations are discussed in \citet{John+15rr}.
We take a different approach and analyze the assumptions necessary to have
performance guarantees for either class of methods. 




Our main contribution is a characterization of the data situations which are
difficult for feature selection algorithms.  This characterization should
provide statistical insight as well as a way to generalize insight gained from
low-dimensional problems. Unfortunately these two are not accomplished in
the same way, which necessitates providing multiple definitions of
approximate submodularity.

\citet{DasK11} introduced a notion of approximate submodularity, measured by
the submodularity ratio, which we will call ``statistical submodularity'' given
its connection to performance guarantees of statistical algorithms. We provide a
characterization of the data situations in which this criteria holds. While the
submodularity ratio is statistically useful, it does not allow insight gained
from low-dimensional problems to be generalized. We provide a stronger
definition for approximate submodularity and show that it yields a lower bound
on the submodularity ratio. In particular, we explain which data conditions
yield approximate submodularity for all feasible two-dimensional regression
problems. While this is restrictive, it yields generalizable bounds and
insights.

As submodularity is a function of model fit, it depends on the response $Y$.
This allows for a broader understanding of problematic correlation structures
and is highly relevant to many simulation settings. From this perspective, not
all deviations from orthogonality are the same. Spectral measures of such
deviations do not always account for this lack of symmetry. Often simulations
are described by their signal to noise ratio without considering the relative
difficulty of different functional forms of the response. Provide an honest
measurement of the difficulty of simulated data case requires considering both
the strength of the signal and the ease with which the signal can be found.

Lastly, we demonstrate that submodularity often appears in statistics
literature, just not by that name. We discuss the restricted eigenvalue
\citep{RasWY10} and conditions for sure independent screening (SIS)
\citep{FanL08}. The discussion highlights the data situations in which the
sparse regression problem (\ref{eqn:sparse-reg}) is solvable by either
approximating the solution or relaxing the problem formulation. Essentially,
achieving an approximate solution to the exact problem is successful is the same
instances in which achieving an exact solution to the approximate problem is
successful. Furthermore, counter-intuitive results from recent conditional
testing literature on forward stepwise and Lasso \citep{Taylor+14} are
explained by deviations from submodularity.

Section \ref{sec:submod} introduces submodularity and our definition of
approximate submodularity. Section \ref{sec:submod-graph} provides a simple
example with only two features to provide intuition about the constraint of
approximate submodularity. Furthermore, it is shown how submodularity can
influence the search path identified by a greedy procedure. We also demonstrate
the effect of signal strength in conjunction with submodularity. If the signal
is strong enough, deviations from submodularity are easier to tolerate because
signal is harder to hide in complex relationships between features. Lastly,
Section \ref{sec:other-assumptions} discusses the connection between
submodularity and more common assumptions in statistics.


\section{Submodularity}
\label{sec:submod}

Submodularity is a condition under which greedy algorithms perform well. In this
section, submodularity is given a statistical interpretation which
begins to reveal its relevance in statistics. We often need to consider a
feature $\X_i$ orthogonal to those currently in the model, $\X_S$. This is
referred to as adjusting $\X_i$ for $\X_S$. The projection operator (hat
matrix), $\Hm_{\X_S} =\Hm_S = \X_S(\X_S^T\X_S)^{-1}\X_S^T$, projects a
vector onto the span of the columns of $\X_S$. Therefore, $\X_i$ adjusted for
$\X_S$ is denoted as residual $\X_{i.S^\perp} = (\I-\Hm_{\X_S})\X_i$. This same
notation holds for sets of features: $\X_A$ adjusted for $\X_S$ is
$\X_{A.S^\perp} = (\I - \Hm_{\X_S})\X_A$.

While assuming \Rs is submodular is uncommon in the statistics literature, an
equivalent formulation has been discussed in the social science literature: the
absence of conditional suppressor variables \citep{DasK08}. It is often
observed features that have positive marginal correlation with the response
can have negative partial correlation in the presence of other features.
Similarly, features can be more significant in the presence
of others than they are in isolation. In these situations, ``suppression'' is
said to have occurred. The words ``suppression'' and ``suppressor variable'' can
be understood through the algebra of adjustment for multiple regression
coefficients.

If $\X$ and $Y$ are standardizes, the coefficient for a feature $\X_i$ in a
simple regression is the correlation between $\X_i$ and $Y$: $r_{Y,i}$. Letting
$C = S\backslash i$ be the other features in the model, the coefficient for
$\X_i$ in a multiple regression is
\begin{IEEEeqnarray*}{rCl}
 \hat\beta_i & = & \frac{\langle Y,\X_{i,C^\perp} \rangle}{\langle
\X_{i,C^\perp}, \X_{i,C^\perp}\rangle}.
\end{IEEEeqnarray*}
Therefore suppression occurs when variability in the feature of interest
that is unrelated to $Y$ is \emph{suppressed} by the other features in the
model.

A suppressor variable is one which, once controlled for, \emph{increases} the
observed significance of another feature. The absence of a conditional
suppressor implies that $\forall S \subset [m]$ and $i,j \notin S$
\begin{IEEEeqnarray*}{rCl}
|\text{Corr}(Y,\X_{i,(S\cup j)^\perp})| & \leq &
|\text{Corr}(Y,\X_{i,(S)^\perp})|.
\end{IEEEeqnarray*}

Suppression is fundamentally the same problem as Simpson's paradox
and Lord's paradox. The distinction arises based on the type of features being
considered. Given features $\X_1$ and $\X_2$, Simpson's paradox
can occur when both features are categorical, Lord's paradox can occur
when one is categorical and the other is numeric, and suppression
can occur when both features are numeric. Any of these paradoxes create
problems with interpreting the influence of features in a regression model.

If one is only interested in prediction, the interpretation of coefficients is
often unimportant. The existence of a suppressor variable does not change the
predictions from a model; however, suppression has significant consequences for
the ability of an algorithm to identify an important feature. In extreme
cases, important features can only be identified as such in the context of
many other features. To extend the simple example given in the introduction,
consider the following set of random variables:
\begin{IEEEeqnarray*}{rCl'rCl}
 \Z & = & N_p(\mathbf{0},\sigma_z \I_p) &
    \mathbf{\epsilon} & = & N_{p-1}(\mathbf{0},\sigma_\epsilon
    \I_{p-1})\numspace \label{eqn:super.start}\\
 X_{1:(p-1)} & = & \Z_{1:(p-1)} + \mathbf{\epsilon} &
    X_p & = & \Z_p - \sum_i^{p-1} \mathbf{\epsilon}_i\\
 Y & = & \sum_{i=1}^{p} X_i = \sum_i^{p} \Z_i \label{eqn:super.end}
\end{IEEEeqnarray*}
Suppose that $\sigma_\epsilon/\sigma_z$ is large enough that the variability in
$\epsilon$ hides any signal that is in $\Z_i$. In this example, any model
with fewer than $p$ features has an \Rs near 0, while using all $p$ features
yields an \Rs of 1. The improvement in fit by adding any single variable is
approximately 0 or 1, depending on which other variables are in the model. This
clearly harms any procedure that solves isolated subproblems. Given the
equivalence between lack of suppression and submodularity, we will use these
term interchangeably. Similarly, subsets of features which violate
Definition \ref{def:submod-1} are instances of supermodularity. Therefore
suppression situations are also supermodular.\footnote{Given that $\mRs\geq0$,
submodular function are also subadditive. Similarly, supermodular ones are
superadditive. While we do not use this terminology, it may be encountered
elsewhere.}
Further implications of the submodularity of \Rs are understood by considering
equivalent definitions of submodular functions. Definition \ref{def:submod-1}
provides the classical definition of submodularity, and two refinements can be
made that merely specify the sets under consideration in increasing detail. For
completeness, all three formulations are provided in Definition
\ref{def:submod-full} and are ordered in terms of specificity.

\begin{defn}[Submodularity]   Let $F: 2^{[m]} \rightarrow \bbR$ be a set
function defined on the the power set of $[m]$. $F$ is submodular iff
\label{def:submod-full}
\begin{enumerate}
 \item (Definition) $\forall A,B \subset [m]$
 \begin{IEEEeqnarray*}{rCl}
  F(A) + F(B) & \geq & F(A \cup B) + F(A \cap B)\\
  \Rightarrow 
  F(A) -F(A \cap B) & \geq & F(A \cup B) -F(B)\\
  \Rightarrow \Delta_{A\cap B}(A) & \geq & \Delta_B(A)
 \end{IEEEeqnarray*}
 \item (First-order difference) $\forall A,B$ such that $A \subset B
\subset [m]$ and $i \in [m]\backslash B$
 \begin{IEEEeqnarray*}{rCl}
  F(A \cup \{i\}) - F(A) & \geq & F(B \cup \{i\}) - F(B)\\
  \Rightarrow \Delta_A(i) & \geq & \Delta_B(i)
 \end{IEEEeqnarray*}
  \item (Second-order difference) $\forall A \subset [m]$ and $i,j \in
[m] \backslash A$
 \begin{IEEEeqnarray*}{rCl}
  F(A \cup \{i\}) - F(A) & \geq & F(A \cup \{i,j\}) - F(A \cup \{j\})\\
  \Rightarrow \Delta_A(i) & \geq & \Delta_{A\cup j}(i)
 \end{IEEEeqnarray*}
 \end{enumerate}
\end{defn}

The definition in terms of first-order differences shows that submodular
functions are similar to concave functions in that they exhibit diminishing
marginal returns. The marginal impact or discrete derivative of adding a
feature to $A$ is larger than that of adding it to $B$ since $A \subset B$. In
terms of optimization, however, they behave like convex
functions and can be efficiently minimized. See \citet{Bach11} for a survey of
this viewpoint. One further simplification is
possible by specifying $B = A \cup \{j\}$, which yields the definition in terms
of second-order differences. This provides the most granular, well-specified
definition of submodularity, and it is the easiest to verify in practice. The
proofs of the equivalence of these definitions are standard and can be found in
many places, for example \citet{Bach11}. Furthermore, when showing the
equivalence of definitions for approximate submodularity, we will be using
proofs of essentially the same form.

In statistical terms, the first- and second-order difference definitions capture
the intuitive notion that correlated features \emph{share} information.
Suppose $\X_S$ is a highly positively correlated set of features where
$\beta_i\geq0,$ $\forall i \in S$. If only $\X_j$, $j\in S$, is added to the
model, it produces a larger marginal improvement in fit than if the entire set
$\X_S$ is included: $\Delta_{\emptyset}(\X_j) \geq \Delta_{S\backslash
j}(\X_j)$. This claim does not hold in general, but does in this example
because it is submodular. Correlation structures which violate this
intuitive notion of shared information are described in Section
\ref{sec:graph-def}.

The above discussion follows from elementary decompositions of simple and
multiple regression coefficients. Let $S = \{i\cup j\}$ and consider the
following models, where subscripts $m$ and $s$ index the model coefficients and
error terms:
\begin{center}
  \begin{tabular}{r|c|c}
  & Multiple Regression & Simple Regression \\\hline
  Model &  $Y = \beta_{0,m}+ \X_i\beta_{i,m} + \X_j\beta_{j,m} +
  \epsilon_s $ & $Y = \beta_{0,s}+ \X_i\beta_{i,s} + \epsilon_m $\\
  Estimated Coefficients & $\hat\beta_{0,m},$ $\hat\beta_{i,m}$ and
  $\hat\beta_{j,m}$ & $\hat\beta_{0,s}$ and $\hat\beta_{i,s}$
  \end{tabular}
\end{center}
The simple regression coefficient can be decomposed into direct and indirect
effects:
\begin{IEEEeqnarray}{rCl}
 \hat\beta_{i,s} & = & \underbrace{\hat\beta_{i,m}}_\text{direct} + 
  \underbrace{\hat\alpha_j\hat\beta_{j,m}}_\text{indirect}.
\label{eqn:coef-decomp3}
\end{IEEEeqnarray}
where $\hat\alpha_j$ is estimated from\[\X_i  =  \alpha_0 + \X_j\alpha_j +
\epsilon.\]
By construction, all terms are positive in equation (\ref{eqn:coef-decomp3}) and
the simple regression coefficient $\hat\beta_{i,s}$ is larger than
$\hat\beta_{i,m}$. Therefore, the marginal impact of adding $\X_i$ is larger in
isolation than in conjunction with $\X_j$. While this is a simplistic example,
it introduces the general insight gained in later sections. In the simplest
case, submodular data requires positively correlated features to have
correlations with the response of the same sign. For example, both must be
negative or positive. Similarly, if features are negatively correlated, their
correlations with the response need to be of opposite sign.


The conditions provided in Definition \ref{def:submod-full} need to be
relaxed in order to capture the continuum of possible scenarios. This will
provide a measure of how signal can ``hide'' in sets of features while not
being visible marginally. This measure is closely connected to assumptions more
commonly discussed in statistics (see Section \ref{sec:other-assumptions}).
There are two conflicting interests when
providing an approximate definition of submodularity. First, it needs to be
statistically meaningful. Such a definition should characterize a relevant
statistical problem that needs to be addressed by many algorithms. Second,
understanding submodularity in spaces with few features should provide
generalizable insight into higher-dimensional problems. Unfortunately, both
goals are not accomplished in the same way. Therefore, two notions of
approximate submodularity are developed and their relationships are described.

Forward stepwise works better if the
influence of a set $S$ can be bounded by the sum of the margin influences of the
elements in it. This can be achieved by applying Definition \ref{def:submod-1}
multiple times to reduce the left hand side to a sum of individual elements. If
$A = \{a_1,\ldots,a_l\} \subset [m]$ and $B = \{b_1,\ldots,b_m\} \subset [m]$,
this yields
\begin{IEEEeqnarray}{rCl}
\sum_{i=1}^l \Delta_{A\cap B}(a_i) + \sum_{i=1}^m
\Delta_{A\cap B}(b_i) & \geq &\Delta_{A\cap B}(A\cup B). \label{eqn:sub-stat}
\end{IEEEeqnarray}
Note that for elements $a_i \in A\cap B$ or $b_i \in A\cap B$
that $\Delta_{A\cap B}(a_i) = \Delta_{A\cap B}(b_i) =0$.

\citet{DasK11} propose a definition of approximate submodularity that
requires equation (\ref{eqn:sub-stat}) to hold approximately by including a
constant $\gamma_{sr}>0$ on the right hand side. This is different than
incorporating the same constant into Definition \ref{def:submod-1} as multiple
applications of the definition are required to produce equation
(\ref{eqn:sub-stat}). For additional simplicity, consider adding the set
$A=\{a_i,\ldots,a_l\} \subset [m]$ to the model $S$. Hence $\Delta_S(a_i)$ is
the
marginal increase in \Rs by adding $a_i$ to model $S$. In this simple case,
$\Delta_S(a_i)$ is the squared partial-correlation between the response $Y$ and
$a_i$ given $S$: $\Delta_S(a_i) = \text{Cor} (Y,a_{i.S}^\perp)^2$.
Therefore, define the vector of partial correlations as $r_{Y,A.S^\perp}
= \text{Cor}(Y, A.S^\perp)$, then the left hand side of
(\ref{eqn:sub-stat}) is $\|r_{Y,A.S^\perp}\|_2^2$. Similarly, if we define
$C_{A.S^\perp}$ as the correlation matrix of $A.S^\perp$ then $\Delta_S(A) =
r_{Y,A.S^\perp}'C_{A.S^\perp}^{-1}r_{Y,A.S^\perp}$.

\begin{defn}
 The submodularity ratio, $\gamma_{sr}$, of \Rs with respect to a set $S$ and $k
\geq 1$ is
 \begin{IEEEeqnarray*}{rCl}
  \gamma_{sr}(S,k) & = & \min_{(T:T\cap S = \emptyset, |T| \leq k)}
\frac{r_{Y,T.S^\perp}'r_{Y,T.S^\perp}}{r_{Y,T.S^\perp}'C_{T.S^\perp}^{-1}r_{
Y,T.S^\perp}}
 \end{IEEEeqnarray*}
\end{defn}
The minimization identifies the worst case set $T$ to add to the model $S$. It
captures how much \Rs can increase by adding $T$ to $S$
(denominator) compared to the combined benefits of adding its elements to $S$
individually (numerator). \Rs is submodular if $\gamma_{sr} \geq 1$ for all $S
\subset [m]$ and $k=2$. Only checking $k=2$ is sufficient due to the
second-order
difference definition of submodularity and is clear from the proofs later in
this section.
\mynotesH{The definition can be made tighter for proofs in terms of how the sets
$L$ and $T$ are defined.}

To not conflate the different notions of approximate submodularity introduced in
this section, $\gamma_{sr}$ will be referred to as the submodularity ratio or
statistical submodularity. It can be used in proofs of the performance of greedy
algorithms \citep{John+15rai,DasK11} and is lower bounded by a sparse eigenvalue
\citep{DasK11}. The connection to spectral quantities is obvious as
$\gamma_{sr}$ is an inverted Rayleigh quotient of the covariance matrix
$C_{T.L^\perp}$. As $C_{T.L^\perp}$ is the Schur complement of $C_{T\cup L}$,
Corollary 2.4 from \citet{Zhang06} proves that $\gamma_{sr}$ is lower bounded by
the minimum eigenvalue of $C_{T\cup L}$. The minimum sparse eigenvalue merely
removes the dependence on the selected sets $L$ and $T$. The connections to
other algorithms that depend on spectral quantities are discussed in Section
\ref{sec:other-assumptions}. 

The submodularity ratio is not appealing from the perspective of submodularity.
It is redefined for different cardinalities $k$ and does not allow information
gained for fixed $k$ to percolate to larger $k$. We now provide a refined
construction of approximate submodularity that produces generalizable
insights. The definitions of approximate submodularity should mirror those of
Definition \ref{def:submod-full}, so that knowledge gained from restrictive,
two-dimensional cases can generalize to higher-dimensional cases. These
equivalent definitions, however, consider submodular functions in a slightly
different context than the submodularity ratio $\gamma_{sr}$. The distinction is
due to bounding the minimum of a set of differences versus the sum of a set of
differences. Clearly bounding the minimum is stronger.

Approximate submodularity is constructed by starting with the second-order
differences definition as it is the most granular and well-specified. Ideally,
the sum of the marginal impact of features considered individually would be
approximately greater than their impact considered jointly. Namely, for some
constant $\gamma > 0$,
\[\Delta_A(i) + \Delta_A(j) \geq \gamma \Delta_A(i,j).\addtag
\label{eqn:submod.sum}\]
This is $\gamma_{A,2}$ after fixing the sets being minimized, but is
unfortunately too weak to generalize to the larger sets considered in Definition
\ref{def:submod-full}. Instead, we must maintain the type of comparisons
considered in the standard definitions.

\begin{defn}[Approximate Submodularity] F is approximately
submodular if there exists constants $\gamma_s$ $\gamma_{s2}$, where
$\gamma_{s2} \geq \gamma_s >0$, such that any of the following hold\\
\begin{enumerate}
 \item (Second order difference) $\forall A \subset [m]$ and $i,j \in
[m] \backslash A$
 \begin{IEEEeqnarray*}{rCl}
  F(A \cup \{i\}) - F(A) & \geq & \gamma_{s2}(F(A \cup \{i,j\}) - F(A \cup
    \{j\}))\\
  \Rightarrow \Delta_A(i) & \geq & \gamma_{s2}\Delta_{A\cup j}(i)
 \end{IEEEeqnarray*}
 \item (First order difference) $\forall A,B$ such that $A \subset B
\subset [m]$ and $i \in [m]\backslash B$
 \begin{IEEEeqnarray*}{rCl}
  F(A \cup \{i\}) - F(A) & \geq & \gamma_s(F(B \cup \{i\}) - F(B))\\
  \Rightarrow \Delta_A(i) & \geq & \gamma_s\Delta_B(i)
 \end{IEEEeqnarray*}
 \item (Definition) $\forall A,B \subset [m]$
 \begin{IEEEeqnarray*}{rCl}
  F(A) -F(A \cap B) & \geq & \gamma_s(F(A \cup B) -F(B))\\
  \Delta_{A\cap B}(A) & \geq & \gamma_s\Delta_B(A)
 \end{IEEEeqnarray*}
\end{enumerate}
\end{defn}

One difference between the definitions for submodularity and approximate
submodularity is that the constant will not be the same in all three cases, as
indicated by our use of $\gamma_s$ and $\gamma_{s2}$; however, if either is
strictly greater than 0, then they both are. We are most interested in
$\gamma_s$, which considers large sets, instead of $\gamma_{s2}$, which only
holds for second order differences. We are able to provide a full account for
$\gamma_{s2}$ though, which yields a conservative lower bound on $\gamma_s$.
Therefore, understanding approximate submodularity in two dimensions gives
generalizable insights. The equivalence of these definitions is proved in the
Appendix.

It is easy to see that $\gamma_{sr} \geq \gamma_{s2}$ in the relevant region in
which $\gamma_{s2}\leq1$. In this region, forward stepwise can perform
poorly. The second-order characterization of $\gamma_{sr}$ in equation
(\ref{eqn:submod.sum}) can be constructed using $\gamma_{s2}$.
\begin{IEEEeqnarray*}{rCl}
 \Delta_A(i) & \geq & \gamma_{s_2}\Delta_{A\cup j}(i)\\
 F(A\cup i) - F(A) & \geq & \gamma_{s2}(F(A \cup \{i,j\}) - F(A\cup\{j\}))\\
 F(A\cup i) - F(A) + \gamma_{s2}(F(A\cup\{j\}) - F(A)) & \geq &
\gamma_{s2}(F(A \cup \{i,j\}) - F(A))\\
 \Rightarrow \Delta_A(i) + \Delta_A(j) & \geq & \gamma_s \Delta_A(i,j).
\end{IEEEeqnarray*}
Where the last line follows since $\gamma_{s2}\leq 1$. The submodularity ratio
fixes the base set; hence the above rearranges the definition of
$\gamma_{s2}$ such that the marginal impact of all features is relative to the
same base set $A$. This yields a bound on the sum of marginal effects, whereas
the $\gamma_{s2}$ is a bound on the minimum of the marginal effects. As
expected, the minimum can yield much worse bounds than the sum; however, as
seen in Section (\ref{sec:graph-def}), not all steps can be taken at this worst
case bound.

\section{Submodularity in 2 Dimensions}
\label{sec:submod-graph}

Attempting to classify types of suppression led \citet{Tze91} to graph
suppression situations that are possible with only two features. These graphs
have unintuitive dimension, double-count data instances, and show impossible
configurations. We analyze the same case, but provide graphs that fully
characterize the set of possible regression problems. This clearly displays the
regions in which $\gamma_{s2}$ and $\gamma_{sr}$ are bounded.

\subsection{Graphing Approximate Submodularity}
\label{sec:graph-def}

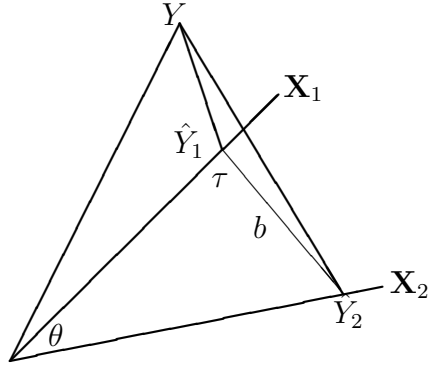
\begin{figure}
\centering
\setlength{\unitlength}{5cm}
\begin{picture}(1,1)
 \thicklines
 \put(0,0){\line(1,2){.447}}
 \put(0,0){\line(1,1){.707}}
 \put(0,0){\line(5,1){.98}}
 \put(.447,.894){\line(1,-3){.112}}
 \put(.447,.894){\line(3,-5){.43}}
 \thinlines
 \put(.559,.559){\line(5,-6){.319}}
 \put(.4,.89){$Y$}
 \put(.72,.7){$\X_1$}
 \put(1,.18){$\X_2$}
 \put(.43,.55){$\hat Y_1$}
 \put(.85,.1){$\hat Y_2$}
 \put(.64,.325){$b$}
 \put(.1,.04){$\theta$}
 \put(.53,.46){$\tau$}
\end{picture}
 \caption{Characterization of possible two-dimensional regression problems: our
data consists of $Y$, $\X_1$, and $\X_2$. $\hat Y_i$ is $Y$ projected
on $\X_i$. The side length from the origin to $\hat Y_i$ is $r_{Yi}$.}
 \label{fig:triangle}
\end{figure}

We parameterize possible regression problems using angles derived from
projecting the response onto individual features. Our data consists of $Y$,
$\X_1$, and $\X_2$ and $\hat Y_i$ as the response $Y$ projected onto $\X_i$. See
Figure \ref{fig:triangle} for an illustration. Since all features have been
normalized, the distance from the origin to $\hat Y_i$ is the correlation
between $Y$ and $\X_i$, $r_{Yi}$. The correlation between explanatory features
is parameterized as $\cos(\theta)$, where $\theta$ is the angle between $\X_1$
and $\X_2$. The relative predictive power of the features is measured by $\tau$,
the angle between $\hat Y_1$ and $\hat Y_2$. Lastly, the strength of the signal
is a function of the length of $b$, the side between $\hat Y_1$ and $\hat Y_2$. 

\begin{figure}[ht]
 \centering
 \includegraphics[width=4in]{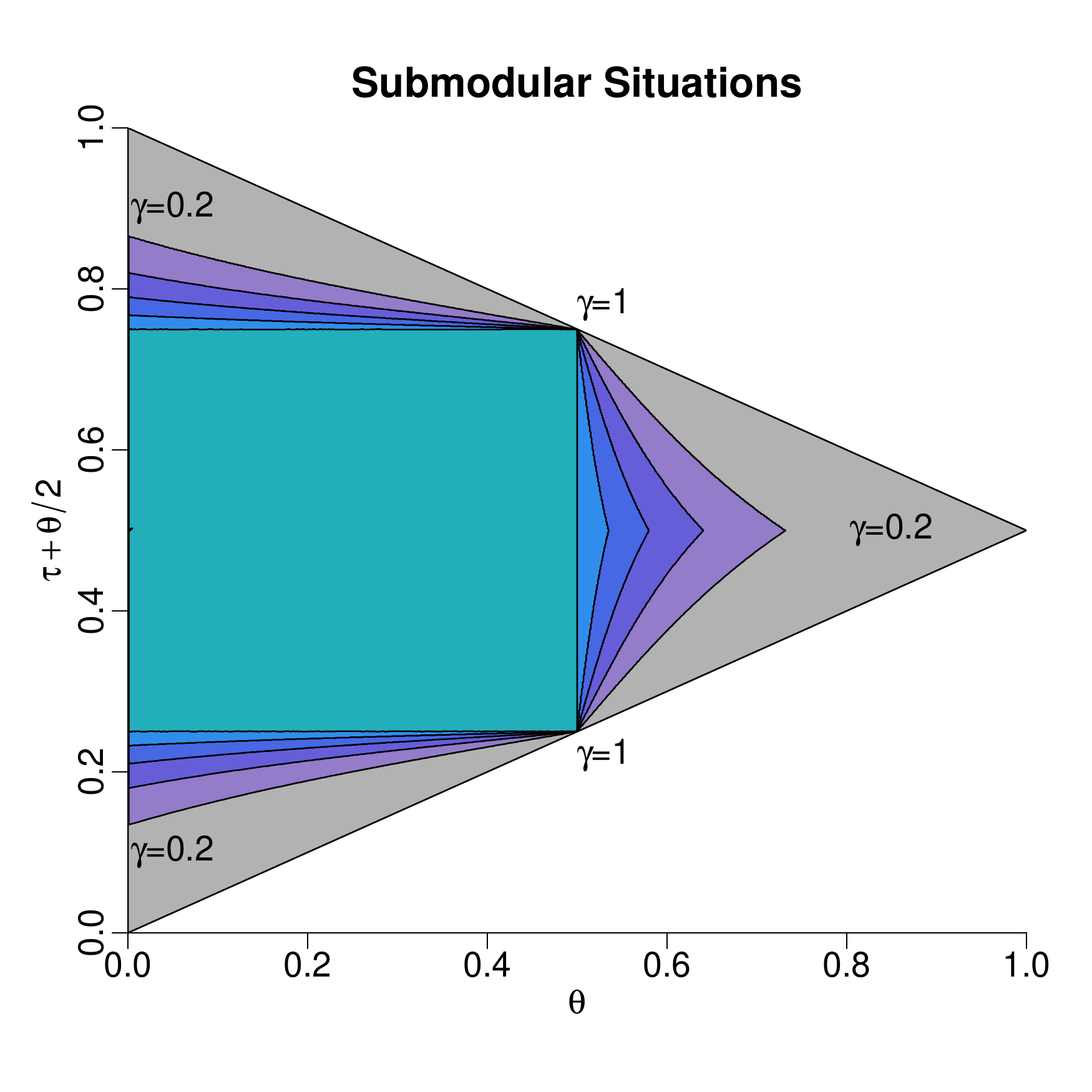}
 \caption{Contour plot of approximate submodularity using second
order differences ($\gamma_{s2}$). Level sets are given for $\gamma_{s2} \in
\{.2,.4,\ldots,1\}$.}
 \label{fig:approx.submod}
\end{figure}

\begin{figure}[ht]
\centering
\includegraphics[width=4in]{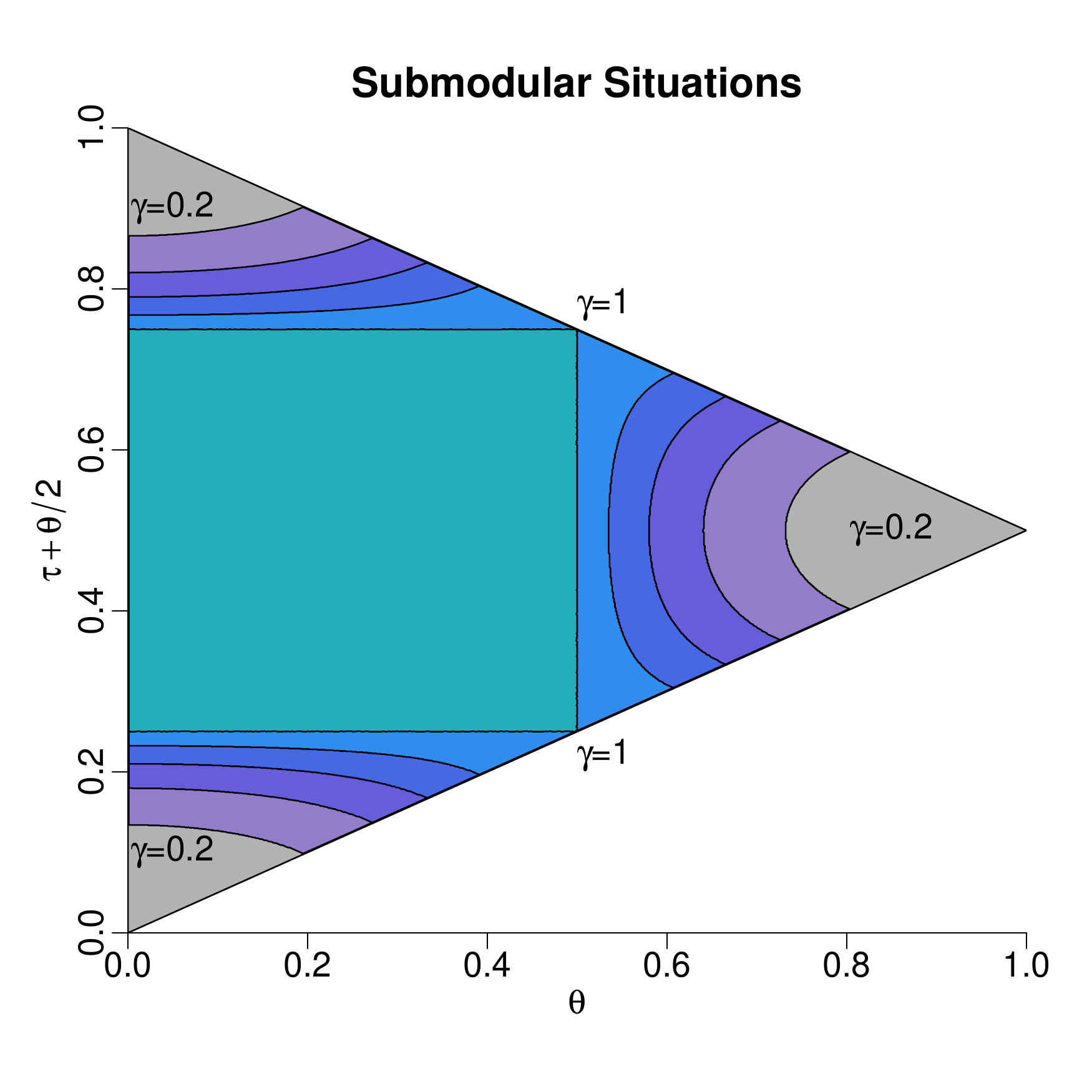}
 \caption{Contour plot of the left hand side of equation (\ref{fig:as.bound}).
The level sets are $\{.2,.4,\ldots,1\}$.}
 \label{fig:as.bound} 
\end{figure}

\begin{figure}[ht]
 \centering
 \includegraphics[width=4in]{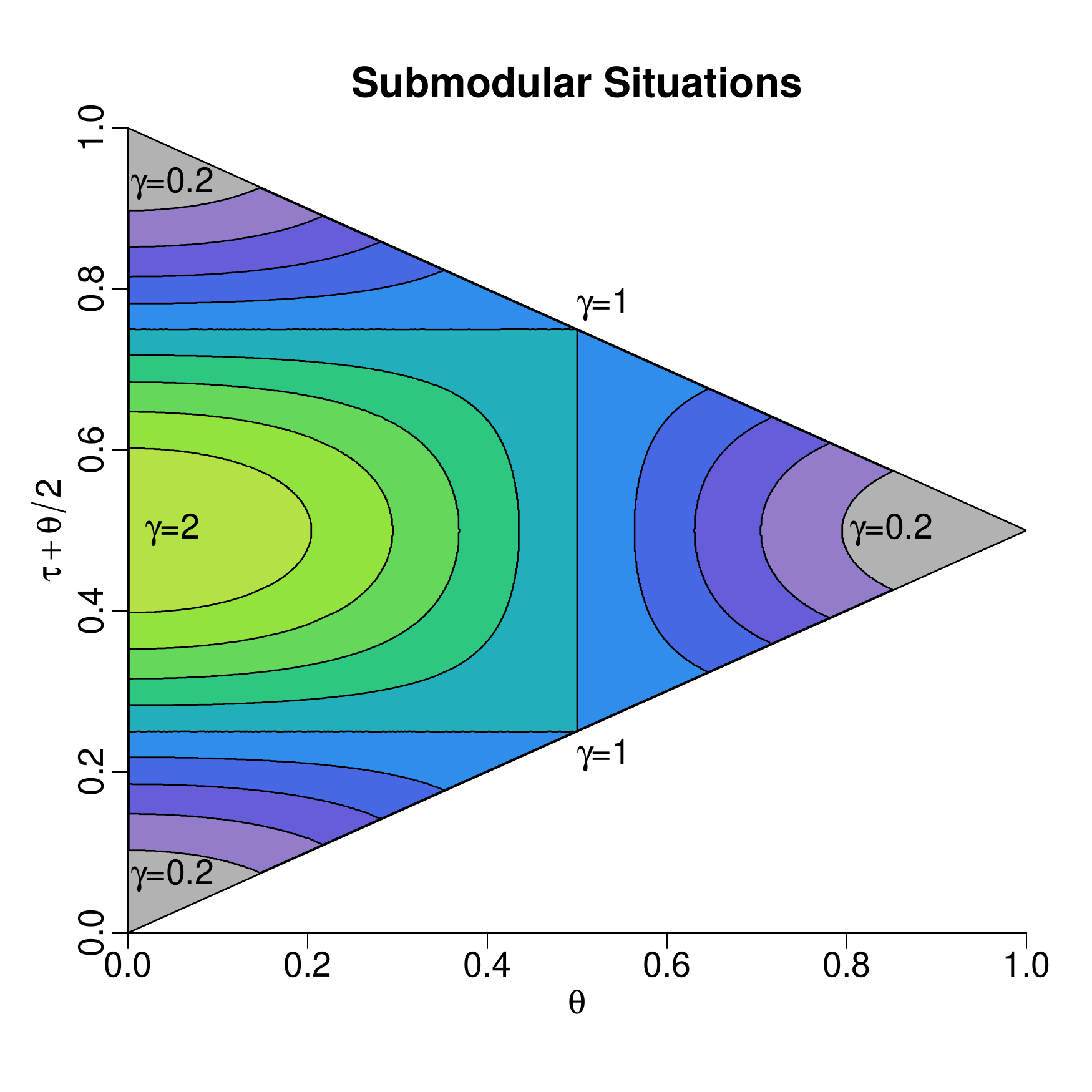}
 \caption{This is a contour plot of the submodularity ratio over the set of
feasible regression problems. Level sets are given for $\gamma_{sr} \in
\{.2,.4,\ldots,2\}$.}
 \label{fig:stat.submod}
\end{figure}

Figures \ref{fig:approx.submod}, \ref{fig:as.bound}, and \ref{fig:stat.submod}
only display $\theta \in [0,\pi]$ and $r_{Yi}\geq0$ because of the symmetries
in submodularity. $\theta > \pi$ is equivalent to $\theta' = (2-\theta)\pi\in
[0,\pi]$ and $r_{Yi} = -r_{Yi}$ for some $i$. The vertical axis has units $(\tau
+\theta/2)\pi$ so that the contour plots are symmetric around $.5\pi$. Figure
\ref{fig:triangle} is an isosceles triangle when $\tau +\theta/2 = .5\pi$,
meaning that both features have the same marginal significance. Therefore,
deviations correspond to one feature being marginally more significant than the
other. Similarly, $\theta=.5\pi$ is the orthogonal case and represents one line
of symmetry on the horizontal axis.

To completely specify the derived triangle in Figure \ref{fig:triangle}, fix a
measure of the signal to noise ratio as this does not represent a
meaningful distinction between models for submodularity. Higher signal just
means that
the effects will be larger. This has the practical impact of being making it
easier to identify a significant effect, but this
discussion is delayed until Section \ref{sec:graph-tstat}. For convenience, we
fix \Rs under the full model: $\mRs_{full} = .5$. All figures are identical
for any value of $\mRs_{full} \in (0,1]$. The length of $b$, the side between
$\hat Y_1$ and $\hat Y_2$, is $\sqrt{(1-r_{12}^2)\mRs_{full}}$.




Figure \ref{fig:approx.submod} is a contour plot of $\gamma_{s2}$ over the set
of feasible regression problems.
It demonstrates that submodularity ($\gamma_{s2} \geq 1$) is only possible when
sign($r_{12}r_{Y1}r_{Y2}$) = 1. This is the intuitive case introduced in Section
\ref{sec:submod}: if features have opposing relationships with the response, we
expect them to be negatively correlated. Since Figure \ref{fig:approx.submod}
displays $r_{Y1}>0$ and $r_{Y2} >0$, submodularity only occurs when the
features are positively correlated. Furthermore, for fixed $r_{12}$ the
maximum $\gamma_{s2}$ occurs when both features have equal marginal effect. As
this is a only a two-feature problem, the joint effects are also equal.
Therefore, the common simulation setting that sets all non-zero coefficients to
the same value maximizes the worse-case step, improving the performance of
feature selection algorithms.

Figure \ref{fig:approx.submod} demonstrates that
while submodularity holds in a large area, relaxing the definition does not
increase the set of problems in a dramatic way; however, this is because
$\gamma_{s2}$ is the single worst case step. Let $\gamma_i$ be
$\Delta(\X_i)/\Delta_{\X_j}(\X_i)$, $i\neq j$, $i,j\in\{1,2\}$, then
$\gamma_{s2}$ is calculated by
\begin{IEEEeqnarray*}{rCl}
\gamma_1 & = & \frac{r_{Y1}^2}{(r_{Y1}^2 - 2r_{Y1}r_{Y2}r_{Y2} +
    r_{Y2}^2r_{12}^2) /(1-r_{12}^2)}\\
\gamma_2 & = & \frac{r_{Y2}^2}{(r_{Y2}^2 - 2r_{Y1}r_{Y2}r_{12} +
    r_{Y1}^2r_{12}^2) /(1-r_{12}^2)}\\
\gamma_{s2} & = & \min(\gamma_1,\gamma_2).
\end{IEEEeqnarray*}
$\gamma_i$ is not symmetric in $\X_1$ and $\X_2$, though given our
interest is in the true model containing both features, it is only important
that one feature appears marginally significant. Importantly, both
features cannot attain the minimum level $\gamma_{s2}$ simultaneously.

To illustrate this, consider bounding the marginal impact of both $\X_1$ and
$\X_2$ using $\gamma_{s2}$. Summing these two inequalities produces
\begin{IEEEeqnarray*}{rCl}
 \frac{\Delta_A(i) + \Delta_A(j)}{\Delta_{A\cup j}(i) + \Delta_{A\cup i}(j)} &
    \geq & \gamma_{s2}\\
 \Rightarrow \frac{\Delta_A(i) + \Delta_A(j)}{2\Delta_A(i,j) -
    \Delta_{A}(i) -\Delta_{A}(j)} & \geq & \gamma_{s2},
\yesnumber\label{eqn:s2-sum}
\end{IEEEeqnarray*}
where the second line just rewrites the first such that the base set is
constant. Figure \ref{fig:as.bound} is a contour plot of the left hand side of
equation (\ref{eqn:s2-sum}). Clearly $\gamma_{s2}$ is a poor bound on this
function, demonstrating that if signal is contained in the joint distribution
of the features, it cannot be hidden from both marginal distributions
simultaneously. It demonstrates that useful properties of submodularity obtain
in much larger region than indicated by $\gamma_{s2}$ due to its
conservativeness.

Lastly, Figure \ref{fig:stat.submod} is a contour plot of the
submodularity ratio $\gamma_{sr}$. It behaves similarly to the bound on the
sum in Figure \ref{fig:as.bound}, though more regularly. There are several
interesting features that can be seen from this graph. First, $\gamma_{sr}$
can
be larger than 1. These are data situations in which forward stepwise achieves
a better bound than the usual $(1-1/e)$ factor off of the optimal. This region
corresponds to cases when the features are highly correlated and have similar
marginal relationships with $Y$. In this case, there is redundancy in our data
and selecting appropriate features is less difficult.

Second, the dependence of $\gamma_{sr}$ on $Y$ is captured by the vertical axis
via $\tau$. Only orthogonal data, $\theta = .5\pi$, is submodular regardless of
$Y$. In this case, the definition of submodularity, equation
(\ref{eqn:submod3}),
holds with equality. This defines a modular function, and it is well known that
the greedy algorithm produces the optimal answer when maximizing a modular
function \citep{Fujishige05}. Due to this dependence on $Y$, $\gamma_{sr}$ is
not
symmetric around the orthogonal case. Obviously the feasible region is not
symmetric, but we consider symmetry in terms of the contours of $\gamma_{sr}$.
The minimum $\gamma_{sr}$ along any vertical strip is achieved at the boundary
of the feasible region. Along the boundaries, submodularity decays at the same
rate when orthogonality is violated with by either positive or negative
correlation. In this way, submodularity is symmetric around the orthogonal case.
This demonstrates the result of \citet{DasK11}, that $\gamma_{sr}$ is lower
bounded by minimum eigenvalues, which occur on the boundary of the feasible
region.
\mynotesH{how does adding noise to decent setup change submodularity?
Assumptions of SIS should maintain this.}

\subsection{Graphing Change in t-Statistics}
\label{sec:graph-tstat}

We now address the issue of \emph{significant} suppression. As the deviation
from submodularity grows, the greedy search path can deviate from the optimal
path; however, slight suppression does not mean that the true model will not be
found. For example, even suppressed features may still be marginally
significant enough to be identified. In this case, the greedy search procedure
has not been harmed.

To analyze these cases, the submodularity ratio can be related to differences in
t-statistics. As in Figure \ref{fig:stat.submod}, consider the contours of the
percentage change in t-statistics caused by different correlation structures.
For clarity, consider the following statistics:
\begin{IEEEeqnarray*}{rCl't}
 \beta_{1m} & = & r_{y1} & m for marginal\\
 \beta_{1j} & = & \frac{r_{y1}-r_{y2}r_{12}}{1-r_{12}^2} & j for joint\\
 t_{1m} & = & \frac{r_{y1}}{\sigma_{im}}\\
 \sigma_{1m}^2 & = & \frac{1-r_{y1}^2}{\sqrt{n-1}}\\
 t_{1j} & = & \frac{(r_{y1}-r_{y2}r_{12})}{(1-r_{12}^2)^{1/2}\sigma_j}\\
 \sigma_j^2 & = & \frac{1}{\sqrt{n-1}}- \frac{r_{y1}^2 -
    2r_{y1}r_{y2}r_{12} + r_{y2}^2}{\sqrt{n-1}(1-r_{12}^2)}
\end{IEEEeqnarray*}
  
\begin{figure}
 \centering
 \includegraphics[width=4in]{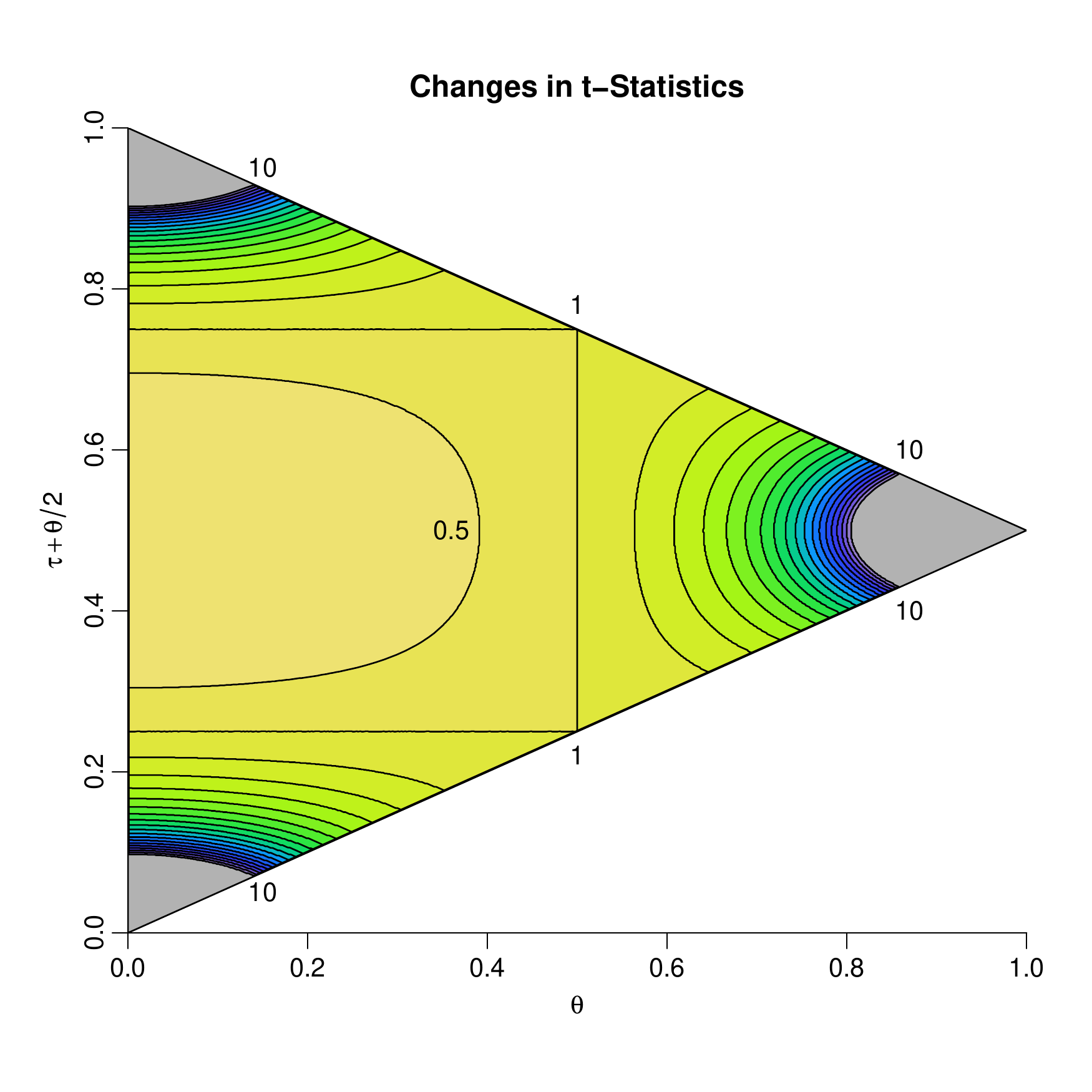}
 \caption{Contour plot of equation (\ref{eqn:tratio}). The contours 
interpolate between .5 and 10 with a step-size of .5.}
 \label{fig:tratio}
\end{figure}

Submodularity requires $t_{1m}^2 \geq t_{1j}^2$. This is a conservative
statement since $\frac{\sigma_{1m}^2}{\sigma_j^2} > 1$. If the features are
jointly highly significant, this becomes very conservative because the ratio is
much larger than 1.
\begin{IEEEeqnarray*}{rCl}
 t_{m1}^2 = \frac{r_{y1}^2}{\sigma_m^2} & \geq & \frac{r_{y1}^2 -
2r_{y1}r_{y2}r_{12} + (r_{y2}r_{12})^2}{(1-r_{12}^2)\sigma_j^2} = t_{j1}^2\\
\Rightarrow r_{y1}^2 & \geq & \frac{r_{y1}^2 - 2r_{y1}r_{y2}r_{12} +
    (r_{y2}r_{12})^2}{1-r_{12}^2}\\
 \Rightarrow r_{y1}^2+r_{y2}^2 & \geq &  \frac{r_{y1}^2 - 2r_{y1}r_{y2}r_{12} +
r_{y2}^2}{1 - r_{12}^2}
\end{IEEEeqnarray*}
Some algebra and incorporating $\gamma_{sr}$ yields the following bound on the
difference between the squared t-statistics:
\begin{IEEEeqnarray*}{rCl}
 \Rightarrow t_{j1}^2 - t_{m1}^2 & \leq & 
    \frac{(1-\gamma_{sr})(r_{y1}^2 - 2r_{12}r_{y1}r_{y2} + r_{y2}^2
)}{1-r_{12}^2}
\end{IEEEeqnarray*}
The previous display ignores the symmetry of the problem: it is not of concern
which of $\X_1$ or $\X_2$ is the suppressed feature, merely that there exists
one. As such, add the corresponding equation for $\X_2$ and divide
by the sum of the marginal t-statistics. This treats $\X_1$ and $\X_2$
symmetrically, and yields
\begin{IEEEeqnarray*}{rCl}
\frac{t_{j1}^2 + t_{j2}^2}{t_{m1}^2+t_{m2}^2} & \leq & 1 +
    \frac{2(1-\gamma_{sr}) (r_{y1}^2 - 2r_{12}r_{y1}r_{y2} + r_{y2}^2
    )}{(1-r_{12}^2)(r_{y1}^2 + r_{y2}^2)} \\
    & = & 2\gamma_{sr}^{-1}-1.\label{eqn:tratio} \IEEEyesnumber
\end{IEEEeqnarray*}
Since equation (\ref{eqn:tratio}) is conveniently written in terms of the
$\gamma_{sr}$, we provide its contour plot in Figure \ref{fig:tratio}. Equation
(\ref{eqn:tratio}) is always positive since $\gamma_{sr} \leq 2$.

The contours of Figure \ref{fig:tratio} are similar to those in
Figure \ref{fig:stat.submod}, but the contours change at different rates.
If $\gamma_{sr} > .8$, then the ratio of squared t-statistics cannot be greater
than 1.5. In this case, if a greedy procedure stops because all
remaining features have a marginal t-statistic less than 2 in absolute value,
neither feature can have a t-statistic larger than 3.46 when considered
jointly. This upper bound is attained when one feature has a joint t-statistic
of 0. If the joint information is split evenly between the two features, the
maximum joint t-statistics are 2.44. Again, it is important to note that \Rs is
not involved in this equation. Therefore submodularity is measuring a
fundamentally different component than the signal-to-noise ratio.
\mynotesH{True in general even though derived only with 2 features?
``sufficient'' construction yields this?}
\mynotesH{Where did \Rs go in this equation? Does it not play a role?}

\section{Connection to Other Assumptions}
\label{sec:other-assumptions}

Some algorithms that leverage assumptions similar to submodularity are the
Lasso, Dantzig selector, and sure independent screening (SIS). It should not be
surprising that the Lasso and forward stepwise are closely connected as the LARS
procedure demonstrates the approximate greedy nature of the Lasso
\citep{Efron+04}.
This similarity extends to the Dantzig selector given that the same assumption
guarantees success of the Lasso and Dantzig selector \citep{BickelRT09}. Lastly,
SIS needs guarantees that information learned from marginal correlations is
sufficient for model selection
\citep{FanL08}. This section describes these procedures, the assumptions used to
demonstrate their success, and their close connection to submodularity.

\subsection{Lasso and Dantzig}

Relaxing the constraint from problem (\ref{eqn:penreg}) from $\|\beta\|_{l_0}$
to $\|\beta\|_{l_1}$ yields the Lasso problem \citep{Tib96}. 
\begin{IEEEeqnarray}{rCl}
\hat\beta_l & = & \text{argmin}_{\beta}\left\{\text{ESS}(\X \beta) +
  \lambda\|\beta\|_{l_1} \right\}
\end{IEEEeqnarray}
This is a convex program and can be efficiently solved using a variety of
algorithms \citep{Efron+04,glmnet}. 

The Dantzig selector \citep{CanT07} optimizes the following linear program
\begin{IEEEeqnarray}{rCl}
\hat\beta_d & = & \text{argmin}_{\beta}\left\{ \|Y-\X \beta\|_\infty +
  \lambda\|\beta\|_{l_1} \right\}
\end{IEEEeqnarray}

Our discussion of these procedures focuses on the assumptions required to
provide bounds on their prediction loss. Many properties have been defined
such as the restricted isometry property \citep{CanT05}, the restricted
eigenvalue constant \citep{BickelRT09,RasWY10}, or the compatibility
condition \citep{Geer07}. For a review of these and related assumptions, see
\citet{GeerB09}. For our purposes, the most important of these is the restricted
eigenvalue, which is defined over a restricted set of vectors that contain
$\hat\beta_l$ and $\hat\beta_d$. Consider a
subset $S \subset \{1,\ldots,p\}$ and constant $\alpha > 1$. Define the set
\begin{IEEEeqnarray*}{rCl}
 C(S;\alpha) & := & \{\beta \in \bbR^p | \|\beta_{S^c}\|_1 \leq \alpha
  \|\beta_S\|_1\}
\end{IEEEeqnarray*}
The restricted eigenvalue of the $p \times p$ sample covariance matrix
$\hat \Sigma = \X^T\X/n$ is defined over $S$ with parameter $\alpha \geq 1$.
\begin{IEEEeqnarray*}{rCl}
 \gamma_{re}^2(\alpha,S) & := & \min \left\{
\frac{\beta'\hat\Sigma\beta}{\|\beta_S\|_2^2} : \beta \in C(S;\alpha)
\right\}
\end{IEEEeqnarray*}
If $\gamma_{re}$ is uniformly lower-bounded for all subsets $S$ with
cardinality $k$, $\hat\Sigma$ satisfies a restricted eigenvalue
condition of order $k$ with parameter $\alpha$. 

The restricted eigenvalue is effectively the submodularity ratio tailored to the
Lasso and Dantzig selectors and generalized to hold for all response vectors
$Y$. Previous work has demonstrated the connection between $\gamma_{sr}$ and
sparse eigenvalues \citep{DasK11}. A sparse eigenvalue with parameter
$k<p$ is
\[\lambda_{\min}(k) = \min_{\delta\in\bbR^k:1\leq \|\delta\|_0 \leq k}
\frac{\delta^T\X^T\X\delta}{n\|\delta\|_2^2}.\]
In order to remove the dependence on $Y$ in the definition of $\gamma_{sr}$,
both the model $S$ of size $k$ and the comparison set $L$ of size $k$ need to
be arbitrary. Therefore, $\gamma_{sr}(S,k) \geq \lambda_{\min}(2k)$. As
discussed in \citet{BickelRT09}, bounding the restricted eigenvalue bounds the
minimum $2k$-sparse eigenvalue. Thus the data conditions under which the Lasso
and Dantzig selector are guaranteed to be successful are stronger than those
under which forward stepwise is. Granted, the form of the guarantees are
significantly different, but of interest is the similarity of the assumptions
required.

The Lasso and Dantzig selector are known to over-estimate the support of $\beta$
\citep{Zou06}, and thus should not be compared to a sparse vector with $k$
non-zero entries. The estimates $\beta_l$ and $\beta_d$ are elements of
$C(S;\alpha)$ with probability close to 1 \citep{BickelRT09}. Therefore, the
bound corresponding to submodularity needs to minimize over $C(S;\alpha)$
instead of truly sparse vectors. Given the looseness of $\lambda_{\min}(2K)$ as
a lower bound on $\gamma_{sr}(S,k)$, we expect a similar looseness exists
between the restricted eigenvalue and the corresponding $Y$-dependent bound.
While it is useful to provide guarantees that do not depend on $Y$, the
potential to produce a better estimate of the crucial constant at runtime may
provide stronger practical performance guarantees. This development could mirror
\citep{BertAR15}.

\subsection{SIS: Sure Independent Screening}

SIS is a correlation learning method in which the marginal correlations between
the response and all features are computed and the features with the largest
$d$
correlations are kept. This can be coupled with subsequent feature
selection algorithms such as SCAD, Dantzig, or Lasso to select a final model
from these $d$ features.
As an additional step, this process can be iterated in much the same way as
stepwise regression: all remaining features are projected off of the selected
set, and the process continues using the residuals from the first model.
Therefore, iterated SIS is similar to a batch greedy method.

\citet{FanL08} split the assumptions for the asymptotic analysis into two
groups: one focuses on parameters of the true regression function and the second
focuses on the sampling distribution of the data. The assumptions on the
true function are stronger than submodularity and the sampling distribution does
not distort this. The most relevant assumption the authors make is the
following:

\begin{assumption}\citet{FanL08}\label{ass:sis}
For some $\kappa$, $0 \leq \kappa < 1/2$, and $c_2$, $c_3 > 0$,
\begin{IEEEeqnarray*}{rCl't'rCl}
 \min_{i \in M_{*}} |\beta_i| & \geq & \frac{c_2}{n^\kappa} & \text{ and }
& 
 \min_{i \in M_{*}} |\text{Cov}(\beta_i^{-1}Y,X_i)| & \geq & c_3.
\end{IEEEeqnarray*}
\end{assumption}


This is of the same form as submodularity by:
\begin{IEEEeqnarray*}{rCl}
  |\text{cov}(\beta_i^{-1}Y,\X_i)| & = & |\beta_i^{-1}||\text{cov}(Y,\X_i)| \\
   & = & |\beta_i^{-1}||r_{Yi}|,
\end{IEEEeqnarray*}
where the second line follows because $\X_i$ and $Y$ are standardized. As
$r_{Yi}$ is the coefficient estimate when $\X_i$ considered marginally,
Assumption \ref{ass:sis} assures that features with non-zero coefficients in
the true model have marginal correlations which are large enough to fall
above the noise level. If $S$ is the true model, this can be
written in a similar form as submodularity as $\Delta(\X_i)\geq
c_3\Delta_{S\backslash i}(\X_i)$. This is the first order difference definition
of submodularity when $A=\emptyset$. Furthermore, this is more restrictive than
statistical submodularity since $\gamma_{sr}>0$ merely requires that there
exists at least one feature which increases the model fit when considered in
isolation. Assumption \ref{ass:sis} requires that all true features increase
model fit when considered in isolation. Therefore, all relevant joint
information is visible from correlations. It is impossible to hide signal in
even two-dimensional subproblems such as those considered in Section
\ref{sec:submod-graph}.

\section{Conclusion}

Submodularity plays an important role in statistics because it characterizes
the difficulty of the \emph{search problem} of feature selection. Assumptions
used to prove the success of the Lasso, such as the restricted eigenvalue and
restricted isometry properties, bound minimum sparse eigenvalues
and hence are stronger assumptions than submodularity. Similarly, SIS requires
true model features to have a bounded discrepancy between
their joint coefficient in the true model and their marginal coefficient from a
simple regression. Bounding this discrepancy is stronger than approximate
submodularity as all true features cannot become vastly \emph{more}
significant in the presence of others. Similarly, worst case data examples can
be crafted by intentionally breaking submodularity. This can be seen in
\citet{Berk+13} and \citet{Miller02}. Due to the importance of submodularity in
discrete optimization, it provides a more theoretically robust assumption than
those more commonly considered in statistics. Furthermore, it characterizes a
different dimension of difficulty than the signal to noise ratio. As such, it is
an important statistic to report in simulated data analyses.

\section{Appendix}

\begin{proof}[Proof of equivalence of Definition \ref{def:submod-full}]
Implications $3. \Rightarrow 2. \Rightarrow 1.$ are clear by
appropriately defining the sets of interest as done when introducing the
definitions of submodularity. To prove the reverse implications, we write
lower-level definitions multiple times using nested sets. Summing these
inequalities and simplifying gives the result.

To prove the first-order definition from the second-order definition, consider
$B = A \cup \{b_1,\ldots,b_k\}$, and apply the second-order definition to sets
$A_i' = A \cup \{b_1,\ldots,b_i\}$. This yields a set of inequalities
\begin{IEEEeqnarray*}{rCl}
 \Delta_A(i) & \geq & \gamma_{s2}\Delta_{A_1'}(i)\\
 \Delta_{A_1'}(i) & \geq & \gamma_{s2}\Delta_{A_2'}(i)\\
 & \vdots & \\
 \Delta_{A_{k-1}'}(i) & \geq & \gamma_{s2}\Delta_{B}(i)\\
 \Rightarrow \Delta_A(i) & \geq & \gamma_{s2}\Delta_{B}(i)
    + (\gamma_{s2}-1)\Delta_{A_1'}(i) + \ldots 
    + (\gamma_{s2}-1)\Delta_{A_{k-1}'}(i)\\
 & \geq & \gamma_{s2}\Delta_{B}(i)
    + \frac{\gamma_{s2}-1}{\gamma_{s2}}\Delta_{A}(i) + \ldots 
    + \frac{\gamma_{s2}-1}{\gamma_{s2}^{k-1}}\Delta_{A}(i) \IEEEyesnumber\\
 & \geq & \left(\gamma_{s2} + (1-\gamma_{s2})
    \frac{1-\gamma_{s2}^{-k}}{1-\gamma_{s2}^{-1}}\right)^{-1} \Delta_B(i)
\end{IEEEeqnarray*}
where the second to last line follows from applying the second order
definition repeatedly to convert $\Delta_{A_i'}$ to $\Delta_{A}$. The constant
in the last line provides a lower bound on $\gamma_s$ and is always strictly
positive if $\gamma_{s2}$ is. It assumes that all of the individual steps are
worst-case steps. As
seen in Section \ref{sec:graph-def}, there are constraints on the number of
steps that can be taken at this worst case level.

Similarly, to prove the standard definition from the first-order definition,
apply the latter multiple times and sum the inequalities to produce $\Delta_A(C)
\geq \gamma_s \Delta_B(C)$. Here $C = \{c_1,\ldots,c_k\}$ and $C \cap A =
\emptyset$. Again, let $A_i' = A \cup \{c_1,\ldots,c_i\}$. Note that since $A
\subset B$ this implies that $B_i' = B \cup \{c_1,\ldots,c_i\}$. This yields a
set of inequalities
\begin{IEEEeqnarray*}{rCl}
  \Delta_A(c_1) & \geq & \gamma_s\Delta_B(c_1)\\
  \Delta_{A_1'}(c_2) & \geq & \gamma_s\Delta_{B_1'}(c_2)\\
  & \vdots &\\
  \Delta_{A_{k-1}'}(c_k) & \geq & \gamma_s\Delta_{B_{k-1}'}(c_k)\\
  \Rightarrow \Delta_{A}(C) & \geq & \gamma_s\Delta_{B}(C)
\end{IEEEeqnarray*}
Where the last line follows by summing the previous lines, canceling most
terms. $\forall S,T \subset [m]$, set $A = S \cap T$, $C = S\backslash T$,
and $B = T$. This yields the result.
\end{proof}


\bibliography{/home/kord/Dropbox/Research/Bib_Stuff/Bibliography}


\end{document}